\documentclass{elsart}
\usepackage[english]{babel}
\usepackage{graphicx}
\usepackage{textcomp}
\usepackage{amsmath}

\usepackage{amssymb}
\usepackage{latexsym}
\newtheorem{proposition}{Proposition}
\newtheorem{remark}{Remark}
\newtheorem{Theorem}{Theorem}
\newtheorem{corollary}{Corollary}
\begin {document}

\begin{frontmatter}

\thanks {Partially supported by RFBR Project Nr. 12-01-00308  and by the government grant of the Russian Federation for
support of research projects implemented by leading scientists,
Lomonosov Moscow State University under the agreement No.
11.G34.31.0054.}

\title
{A representation of solutions to a scalar conservation law in
several dimensions}




\author[1]{Sergio Albeverio }
\author[2]{Olga Rozanova\corauthref{cor1}}
\ead{rozanova@mech.math.msu.su} \corauth[cor1]{Corresponding author}

\address[1]{Universit\"{a}t Bonn,
Institut f\"{u}r Angewandte Mathematik, Abteilung f\"{u}r
Stochastik, Endenicher Allee 60, D-53115 Bonn; HCM, SFB611 and IZKS,
Bonn; BiBoS, Bielefeld--Bonn; CERFIM, Locarno}
\address[2]{Department of Mechanics and Mathematrics, Moscow
State University, Moscow 119991 Russia}

\date{\today}

\begin{abstract}
We find a representation of smooth solutions to the Cauchy problem
for a scalar multidimensional conservation law as small diffusion
limit of a stochastic perturbation along characteristics. It helps,
in particular, to study the process of singularities formation.
Further, we introduce an associated system of balance laws that can
be interpreted as describing the motion of a continuum with some
specific pressure term. This term arises only after the instant when
the solution to the initial Cauchy problem looses its smoothness.
Before this instant the system coincides partly with the one known
as pressure free gas dynamics.

\end{abstract}

\begin{keyword}scalar conservation law \sep the Cauchy problem \sep
representation of solution  \sep associated conservation laws  \sep
stochastic perturbation

\MSC 35L65 \sep 35L67
\end{keyword}

\end{frontmatter}




                            %


\section*{Introduction}

We consider the initial value problem
\begin{equation}\label{burgers}
 u_t + \sum\limits_{i=1}^n a_i(t,x,u) \,
u_{x_i}=0,\quad u(0,x)=u_0(x),\quad u_0(x)\in C_b^1({\mathbb R}^n;
{\mathbb R}),\end{equation} where $\, t\in \mathbb R_+,\,x\in
\mathbb R^n,$
 $a_i(t,x,u), \,i=1,...,n,$ is a real-valued $C^1$ function defined on
some open subset of $({\mathbb R}_+ \times {\mathbb R}^n \times
{\mathbb R})$. For a technical reason the functions $a_i(t,x,u)$ are
assumed to grow at infinity not quicker than a linear function in
$u$.

A particular important special case is given by  the scalar
conservation law in the form
$$
u_t+{\rm div} F(t,u)= 0,
$$
where $F(t,\cdot)=(F_1(t,\cdot),...,F_n(t,\cdot))$ is a $C^2$
vector-function defined on some open subset of $ {\mathbb R}$, for
any $t\in {\mathbb R}_+$, and $a_i(t,u)=\frac{\partial
F_i(t,u)}{\partial u},\, i=1,...,n.$

The main aim of this paper is to obtain an asymptotic  formula for
the solution of the Cauchy problem (\ref{burgers}) for the case of a
scalar conservation law. The formula is obtained by the limit for
vanishing perturbation of the corresponding stochastically modified
equation (small diffusion limit).

Nevertheless, let us first consider the general case.

Let us write the associated characteristic ODE:
$$
\frac{dx_i}{dt} = a_i(t,x,u),\quad \frac{du}{dt}=0, \quad i=1,...,n.
$$
Its stochastic analog is
\begin{equation}
\label{SDU}dX_i(t)=a_i(t,X(t), U(t))dt+\sigma_1 d(W_i^1)_t,\quad
dU(t)=\sigma_2 d(W^2)_t,
\end{equation}
$$X_i(0)=x_i,\quad U(0)=u,\quad t>0,$$
$i=1,...,n,$\,  $X(t)$ and $U(t)$ are considered as random variables
with given initial distributions, $(X(t),U(t))$ runs in the phase
space $\mathbb{R}^n\times\mathbb{R}^1,$ $\sigma_1$ and $\sigma_2$
are nonnegative constants such that $|\sigma|\ne 0$
($\sigma=(\sigma_1,\sigma_2)$) and
$((W^1)_t,(W^2)_t)=(W^1_1,\dots,W^1_n, W^2)_t $ is an $n+1$ -
dimensional Brownian motion, i.e. the $W^1_i, W^2$, $i=1,...,n$, are
independent one-dimensional standard Brownian motions.

Let  $P(t,dx,du),\, t\in \mathbb R_+,\,x\in \mathbb R^n,$ be  the
probability  of the joint distribution of the random variables
$(X,U), $ subject to the initial data
\begin{equation}
\label{P0} P_0(dx,du)=\delta_u(u_0(x))\,\rho_0(x) dx,
\end{equation}
where $\rho_0$ is a bounded nonnegative function from $C ({\mathbb
R}^n)$ and $dx$ is Lebesgue measure on ${\mathbb R}^n$,
$\delta_u$ is Dirac measure concentrated on $u$.  $\,P(t,dx,du)$ has
the form $P(t,x,du)dx$, where $P(t,x,du)$ is a positive measure with
respect to $u$ and a function with respect to $x$ (density function
of $P(t,dx,du)$ with respect to Lebesgue measure).

We look at  $P=P(t,dx,du)$ as a generalized function (distribution)
with respect to the variable $u$. It satisfies the Fokker-Planck
equation
\begin{equation}
\label{Fok-Plank} \dfrac{\partial P}{\partial
t}=\left[-\sum\limits_{k=1}^n \dfrac{\partial}{\partial
x_k}\,a_k(t,x,u)
 +\sum\limits_{k=1}^n\dfrac12\sigma_1^2\dfrac{\partial^2}{\partial
x_k^2}+\sum\limits_{k=1}^n\dfrac12\sigma_2^2\dfrac{\partial^2}{\partial
u_k^2}\right] P,
\end{equation}
subject to the initial data \eqref{P0}.

There is a standard procedure for finding the fundamental solution
for \eqref{Fok-Plank} (see, e.g. \cite{Friedman}). This procedure
consists in a  reduction of the equation to a Fredholm integral
equation, the solution of which can be found in the form of series.
We are going to show that for  $a(t,x,u)=a(t,u)$ one can also find
an explicit solution to the Cauchy problem \eqref{Fok-Plank},
\eqref{P0}.

Let us introduce, still in the general case, the functions, for
$t\in \mathbb R_+,\,x\in \mathbb R^n$, depending on
$\sigma=(\sigma_1, \sigma_2)$:
\begin{equation}
\label{rho_s} \rho_\sigma(t,x) = \int\limits_{\mathbb{R}}P(t,x,du),
\end{equation}
\begin{equation}
\label{u_s} \quad
u_{\sigma}(t,x)=\dfrac{\int\limits_{\mathbb{R}}uP(t,x,du)}{\int\limits_{\mathbb{R}}P(t,x,du)},
\end{equation}
\begin{equation}
\label{a_s} \quad
a_{\sigma}(t,x)=\dfrac{\int\limits_{\mathbb{R}}a(t,x,u)
P(t,x,du)}{\int\limits_{\mathbb{R}}P(t,x,du)},
\end{equation}
the integrals in the numerator being assumed to exist in the
Lebesgue sense.

It will readily be observed that $u_{\sigma}(0,x)=u_0(x)$ and
$a_{\sigma}(0,x)=a(0,x,u_0(x))$.

We denote
$$\bar \rho (t,x)=\lim\limits_{\sigma\to 0} \rho_\sigma(t,x),\quad
\bar u(t,x)=\lim\limits_{\sigma\to 0} u_\sigma(t,x),\quad \bar
a(t,x)=\lim\limits_{\sigma\to 0} a_\sigma(t,x),$$ provided these
limits exist.


\section{Case of a conservation law}

Now we dwell on the simpler case of a conservation law, where
$a=a(t,u)$. Here the equation \eqref{Fok-Plank} can be solved
explicitly. Moreover, for the sake of simplicity we set $\sigma_2=0$
and denote $\sigma_1=\sigma$.

\begin{proposition}\label{Prop} If $a=a(t,u)$, then problem
\eqref{Fok-Plank}, \eqref{P0} has the following solution:
\begin{equation}\label{s_plotn}
 P(t,x,du) =\dfrac1{(\sqrt{2\pi
t}\sigma)^n}\int\limits_{\mathbb{R}^n}\,\delta_u(u_0(y))\,\rho_0(y)\,
e^{-\frac{\sum\limits_{i=1}^{n}|\int\limits_0^t{a_i(\tau,u_0(y))d\tau+y_i-x_i|^2}}{2\sigma^2t}}
\,dy,
\end{equation}
$\quad t\ge 0, \,x\in \mathbb{R}^n,$ or, in other words,
\begin{equation}\label{s_plotn_1}
 \int\limits_{\mathbb R} \, \phi(u)\,P(t,x,du) =\dfrac1{(\sqrt{2\pi
t}\sigma)^n}\int\limits_{\mathbb{R}^n}\,\phi(u_0(y))\,\rho_0(y)\,
e^{-\frac{\sum\limits_{i=1}^{n}|\int\limits_0^t{a_i(\tau,u_0(y))d\tau+y_i-x_i|^2}}{2\sigma^2t}}
\,dy,
\end{equation}
for all $\phi(u)\in C_0({\mathbb R})$.

\end{proposition}

{\sc Proof}. We act as in \cite{Roz2}, \cite{Roz3}. Namely, we apply
the Fourier transform to $P(t,x,du)$ in (\ref{Fok-Plank}),
(\ref{P0}) with respect to the  variable $x$    and obtain the
Cauchy problem for the Fourier transform
$\tilde{P}=\tilde{P}(t,\lambda,du)$ of $P(t,x,du)$:
\begin{equation}
\label{preobr_Fok-Plank} \dfrac{\partial \tilde{P}}{\partial
t}=-(\dfrac12\sigma^2|\lambda|^2+i(\lambda, a(t,u)))\tilde{P},
\end{equation}
\begin{equation}
\label{preobr_P0}\tilde{P}(0,\lambda,du)=\int\limits_{\mathbb{R}^n}e^{-i(\lambda,y)}\delta_u(u_0(y))\rho_0(y)dy,
\qquad \lambda \in {\mathbb R}^n.
\end{equation}
Equation (\ref{preobr_Fok-Plank}) can easily be  integrated and we
obtain the solution given by the following formula:
\begin{equation}
\label{preobr_P}\tilde{P}(t,\lambda,du)=\tilde{P}(0,\lambda,du)e^{-\frac12\sigma^2|\lambda|^2t
+i \int\limits_{0}^t (\lambda,a(\tau,u)) d \tau}.
\end{equation}
The inverse Fourier transform (in the distributional sense) allows
to find the density function $P(t,x,du),\,t>0$:
$$P(t,x,du)=\dfrac1{(2\pi)^{n}}
\int\limits_{\mathbb{R}^n}e^{i(\lambda,x)}\tilde{P}(t,\lambda,du)\,d\lambda
=$$
$$=\dfrac1{(2\pi)^{2n}}\int\limits_{\mathbb{R}^n}e^{i(\lambda,x)}
\left(\int\limits_{\mathbb{R}^n}e^{-i(\lambda,y)}e^{-i
\int\limits_{0}^t (\lambda,a(\tau,u)) d
\tau}\,\delta_u(u_0(y))\,\rho_0(y)dy\right)
\,e^{-\frac12\sigma^2|\lambda|^2t}d\lambda=$$
$$=\dfrac1{(2\pi)^{n}}\int\limits_{\mathbb{R}^n}\delta_u(u_0(y))\,\rho_0(y)
\int\limits_{\mathbb{R}^n}e^{-\frac12\sigma^2t\left(\lambda-\frac{i|x-\int\limits_{0}^t
a(\tau,u) d \tau-y|}{\sigma^2t}\right)^2-\frac{|\int\limits_{0}^t
a(\tau,u) d \tau +y-x|^2}{2\sigma^2t}}d\lambda dy=$$
$$
=\dfrac1{(\sqrt{2\pi
t}\sigma)^n}\int\limits_{\mathbb{R}^n}\,\delta_u(u_0(y))\,\rho_0(y)\,e^{-\frac{|\int\limits_{0}^t
a(\tau,u_0(y)) d \tau+y-x|^2}{2\sigma^2t}}dy,\quad t\ge
0,\,x\in \mathbb{R}^n. $$
The third equality is satisfied by Fubini's theorem, which can be
applied by the absolute integrability and the bound on the function
involved. Thus, the proposition is proved.

\begin{remark} In the general case $\sigma_2\ne 0$ an analogous
formula can be obtained in a similar way.
\end{remark}

\begin{corollary} The functions ${\rho}_\sigma$, ${u}_\sigma$ and
$a_\sigma$ defined in \eqref{rho_s} -- \eqref{a_s}  can be
represented by the following formulae:
\begin{equation}
\label{rho_s_rep} {\rho}_\sigma(t,x)=
{\int\limits_{\mathbb{R}^n}\rho_0(y)\,e^{-\frac{\sum\limits_{i=1}^{n}
|\int\limits_0^t{a_i(\tau,u_0(y))d\tau+y_i-x_i|^2}}{2\sigma^2t}}
dy},
\end{equation}
\begin{equation}
\label{u_s_rep} {u}_\sigma(t,x)=
\dfrac{\int\limits_{\mathbb{R}^n}u_0(y)\rho_0(y)\,e^{-\frac{\sum\limits_{i=1}^{n}|\int\limits_0^t{a_i(\tau,u_0(y))d\tau+
y_i-x_i|^2}}{2\sigma^2t}} dy}
{\int\limits_{\mathbb{R}^n}\rho_0(y)\,e^{-\frac{\sum\limits_{i=1}^{n}
|\int\limits_0^t{a_i(\tau,u_0(y))d\tau+y_i-x_i|^2}}{2\sigma^2t}}
dy},
\end{equation}
\begin{equation}
\label{a_s_rep} {a}_\sigma(t,x)=
\dfrac{\int\limits_{\mathbb{R}^n}a(t,u_0(y))\rho_0(y)\,
e^{-\frac{\sum\limits_{i=1}^{n}|\int\limits_0^t{a_i(\tau,u_0(y))d\tau+y_i-x_i|^2}}{2\sigma^2t}}
ds}
{\int\limits_{\mathbb{R}^n}\rho_0(y)\,e^{-\frac{\sum\limits_{i=1}^{n}
|\int\limits_0^t{a_i(\tau,u_0(y))d\tau+y_i-x_i|^2}}{2\sigma^2t}}
ds}.
\end{equation}
\end{corollary}

{\sc Proof.}  The result is obtained by substitution of $P(t,x,du)$
as given by \eqref{s_plotn} in \eqref{rho_s}, \eqref{u_s} and
\eqref{a_s}.

\subsection{Asymptotic formula for smooth solutions}

Let us define the following subset $\Lambda$ of ${\mathbb R}$:
\begin{equation}\label{Lambda}
t\in \Lambda \quad  \mbox{if } \quad \inf\limits_{y\in
\mathbb{R}^n}\int\limits_0^t{\sum\limits_{i=1}^{n}
(a_i)_u(\tau,u_0(y))(u_0(y))_{y_i}}d\tau>-1,
\end{equation}
where $u_0\in C^1_b({\mathbb R}^n)$. It is not difficult to show the
$\Lambda$ is an open set.
 Denote $t_*(u_0)=\sup \Lambda$.

The following  theorem holds:

\begin{Theorem}\label{T1} Let $u(t,x)$ be a solution to the Cauchy problem
\begin{equation}\label{burgers_cl}
 u_t + \sum\limits_{i=1}^n a_i(t,u) \,
u_{x_i}=0,\quad u(0,x)=u_0(x),
\end{equation}
where $a_i, \,i=1,...,n,$ are $C^1$- functions defined on some open
subset of $({\mathbb R}_+ \times {\mathbb R})$ and $u_0\in
C^1_b({\mathbb R}^n)$. Assume $t_*(u_0)=\sup \Lambda>0$, $\Lambda$
being defined by  (\ref{Lambda}). Then for $t\in [0, t_*(u_0)),$
$$u(t,x)=\bar u(t,x)=\lim\limits_{\sigma\to 0} u_\sigma(t,x),$$
where $u_\sigma(t,x)$ is given by (\ref{u_s}) and the limit exists
pointwise.
\end{Theorem}

{\sc Proof}. The proof is similar to the one given in \cite{Roz3}
for a related problem. According to the classical theory (see, e.g.
\cite{Dafermos}), Theorem 5.1.1), the solution $u$ of
\eqref{burgers_cl} exists on some maximal interval $[0, T), \,
T\le\infty$ and is a $C^1$ - smooth function. Since $u$ is constant
along characteristics, its value at any point $(t, x)$, with $x\in
{\mathbb R}^n, $  $t\in {\mathbb R}_+, $ satisfies the implicit
relation
\begin{equation}\label{implicit}
 u(x, t) = u_0(x - \int\limits_0^t\, a(\tau,u) d\tau ).
 \end{equation}
In particular, the range of $u$ coincides with the range of $u_0$.

Differentiating \eqref{implicit} yields
\begin{equation}\label{deriv}
\partial_{x_i} u(t,x)=\frac{\partial_{y_i} u_0(y)}{1+\int\limits_0^t{\sum\limits_{i=1}^{n}
(a_i)_u(\tau,u_0(y))(u_0(y))_{y_i}}d\tau},\quad y=x -
\int\limits_0^t\, a(\tau,u) d\tau.
\end{equation}
This imply $T=t_*(u_0)$. If $0<t_*(u_0)<+\infty$, then the solution
to the Cauchy problem blows up at the instant $t_*(u_0)$. Otherwise,
the solution keeps its smoothness for all $t>0$.


The formula \eqref{u_s}  implies, using the weak convergence of
measures and the fact that $\rho_0$ and $u_0$ are continuous and
bounded and independent of $\sigma$
$$\lim\limits_{\sigma\rightarrow
0}{u}_\sigma(t,x)\,=\,\dfrac{\int\limits_{\mathbb{R}^n}u_0(y)\rho_0(y)\lim\limits_{\sigma\rightarrow
0}\frac1{(\sqrt{2\pi t}\sigma)^n}e^{-\frac{|\int\limits_0^t
a(\tau,u_0(y)) d\tau
+y-x|^2}{2\sigma^2t}}dy}{\int\limits_{\mathbb{R}^n}\rho_0(y)\lim\limits_{\sigma\rightarrow
0}\frac1{(\sqrt{2\pi t}\sigma)^n}e^{-\frac{|\int\limits_0^t
a(\tau,u_0(y)) d\tau +y-x|^2}{2\sigma^2t}}dy}=$$
$$
\,\dfrac{\int\limits_{\mathbb{R}^n}u_0(y)\rho_0(y)\delta_{p(t,x,y)}dy}{\int\limits_{\mathbb{R}^n}\rho_0(y)
\delta_{p(t,x,y)}dy},
$$
with \begin{equation}\label{p}
p(t,x,y): = \int\limits_0^t
a(\tau,u_0(y)) d\tau +y-x,
\end{equation}
 where $\delta_p$ is the Dirac measure at $p\in {\mathbb
R}^n.$ We can use locally the implicit function theorem and find
$y=y_{t,x}(p)$ from $p(t,x,y)$. The condition for existence of this
function is the invertibility of the matrix
$$
C_{ij}(t,y)=\frac{\partial p_i(t,x,y)}{\partial y_j}, \,
i,j=1,...,n.
$$
This matrix fails to be invertible for $t=t_*(u_0).$
For $t<t_*(u_0)$
$$\bar u(t,x)=\lim\limits_{\sigma\rightarrow
0}{u}_\sigma(t,x)\,=\,$$$$
\dfrac{\int\limits_{\mathbb{R}^n}u_0(y_{t,x}(p))\rho_0(y_{t,x}(p))\,\det
(C(t,y_{t,x}(p)))^{-1}\,\delta_p\,(dy_{t,x})}
{\int\limits_{\mathbb{R}^n}\rho_0(y_{t,x}(p))\,\det(C(t,y_{t,x}(p)))^{-1}\,\delta_p\,(dy_{t,x})\,}=\,u_0(y_{t,x}(0)).
$$
Let us introduce the new notation  $y_0(t,x)\equiv y_{t,x}(0).$ Then
\eqref{p} implies the following vectorial equation:
\begin{equation}
\label{usl}\int\limits_0^t a(\tau,u_0(y_0(\tau,x))) d\tau
+y_0(t,x)-x=0, \quad t\ge 0,\quad x\in {\mathbb R}^n.
\end{equation}
Let us show that $u(t,x)=u_0(y_0(t,x))$ satisfies equation
\eqref{burgers}, that is
\begin{equation}\label{Burgsubs}\sum\limits_{j=1}^n\,\partial_j
(u_{0})(y_{0,j})_t\,+\,\sum\limits_{j,k=1}^n a_j(t,u_{0})
\partial_k(u_{0})(y_{0,k})_{x_j}=0.
\end{equation}
and $u_0(y_0(0,x))=u_0(x)$. Here we denote by  $y_{0,i}$ the $i$ -
th components of the vector  $y_0$.

For $t<t_*(u_0)$ we can differentiate (\ref{usl}) with respect to
$t$ and $x_j$ to get the matrix equations:
$$\sum\limits_{j=1}^n\,C_{ij}\,(y_{0,j})_t +u_{0,i}=0,\, \quad
i=1,...,n,$$ and
$$\sum\limits_{k=1}^n\,C_{ik}\,(y_{0,k})_{x_j}+\delta_{ij}=0,
\, \quad i,j=1,...,n,$$ where $\delta_{ij}$ is the Kronecker symbol.
The equations imply
\begin{equation}\label{subst}
(y_{0,j})_t\,=\,-\,\sum\limits_{i=1}^n\,(C^{-1})_{ij}\,u_{0,i},\qquad
(y_{0,k})_{x_j}\,=\,-\,(C^{-1})_{jk}.
\end{equation}
 It remains now only to
substitute (\ref{subst}) into (\ref{Burgsubs}) to see that $\bar
u(t,x)$ satisfies the first equation in \eqref{burgers_cl}.

Further,  (\ref{usl}) implies $u_0(y_0(0,x))=u_0(x)$, thus  Theorem
\ref{T1}  is proved.

\begin{remark}
For $a_i=a_i(u)$ (i.e. $a_i$ is independent of variable $t$) we have
the Conway' criterium \cite{Conway}:
$$
t_*(u_0)= \sup\limits_{y\in
\mathbb{R}^n}\left(-\frac{1}{\sum\limits_{i=1}^{n}
(a_i)_u(u_0(y))(u_0(y))_{y_i}}\right).
$$
Note that if the denominator vanishes, then $t_*(u_0)=\infty$ and
the solution does not blow up. If $t_*(u_0)<0$, then the solution is
 globally smooth for $t\ge 0$, as well.
\end{remark}

\begin{proposition} \label{Prop2} Under the assumptions of Theorem \ref{T1} the vector
$$\bar a(t,x)=\lim\limits_{\sigma\to 0}
a_\sigma(t,x),$$ where $a_\sigma(t,x)$ is given by (\ref{a_s}),
solves the multidimensional Burgers equation
$$
(\bar a(t,u))_t + (\bar a(t,u), \nabla) \bar a(t,u)=0
$$
with initial data $\bar a(0,u(0,x))=a(0,u_0(x)).$
\end{proposition}

{\sc Proof.} This fact follows directly from Proposition 2.1 of
\cite{Roz3}.

\begin{remark} The introduction of a small perturbation of
deterministic equation to study then the original equation in the
limit of vanishing noise has appeared in several contexts,
particularly for equations of the reaction-diffusion type, see, e.g.
\cite{ADM}, \cite{Fl} and references therein.

\end{remark}

\subsection{Associated system of balance laws}

Now we consider the following  question: what system of equations do
the triple $(\rho_\sigma, u_\sigma, a_\sigma)$ and its limit
$(\bar\rho, \bar u, \bar a)$ satisfy before and after the blow up
time $t_*(u_0)$?

The following proposition holds:

\begin{proposition} \label{prop3} The functions $\rho_\sigma$,  $u_\sigma$ and $a_\sigma$, given by \eqref{rho_s}, \eqref{u_s}
and \eqref{a_s}, satisfy for $t\ge 0$ the following PDE system:
\begin{equation}
\label{sist_obw1}\dfrac{\partial\rho_\sigma}{\partial t}\,+\,{\rm
div}_x
(\rho_{\sigma}a_{\sigma})\,=\,\dfrac12\sigma^2\sum\limits_{k=1}^{n}\dfrac{\partial^2\rho_\sigma}{\partial
x_k^2},
\end{equation}
\begin{equation}
\label{sist_obw_u}\dfrac{\partial(\rho_{\sigma}u_{\sigma})}{\partial
t}\,+\, {\rm div}_x(\rho_{\sigma}\,
u_{\sigma}\,a_{\sigma})\,=\,\dfrac12\sigma^2\sum\limits_{k=1}^{n}
\dfrac{\partial^2(\rho_{\sigma}u_{\sigma})}{\partial
x_k^2}\,-I^u_{\sigma },
\end{equation}
where
$$I^u_{\sigma}=\int\limits_{\mathbb{R}^n}(u-u_{\sigma}(t,x))((a(t,u)-a_{\sigma}(t,x)),\nabla_x
P(t,x,du));$$
\begin{equation}
\label{sist_obw2}\dfrac{\partial(\rho_{\sigma}a_{\sigma,i})}{\partial
t}\,+\, {\rm
div}_x(\rho_{\sigma}a_{\sigma,i}\,a_{\sigma})\,=\,\dfrac12\sigma^2\sum\limits_{k=1}^{n}
\dfrac{\partial^2(\rho_{\sigma}a_{\sigma,i})}{\partial
x_k^2}\,-I^a_{\sigma,i },
\end{equation}
 $i=1,..,n,\,$ where
$$I^a_{\sigma,i}=\int\limits_{\mathbb{R}^n}(a_{i}(t,u)-a_{\sigma,i}(t,x)((a(t,u)-a_{\sigma}(t,x)),\nabla_x
P(t,x,du))+$$
$$\int\limits_{\mathbb{R}^n}a_{i}(t,u))_t
P(t,x,du).$$

\end{proposition}

{\sc Proof.} The equation \eqref{sist_obw1} follows from the
Fokker-Planck equation (\ref{Fok-Plank}) directly.

Let us prove \eqref{sist_obw2} (the derivation of \eqref{sist_obw_u}
is analogous).  We note that the definitions of ${a}_\sigma(t,x)$
and $\rho_\sigma(t,x)$ imply
\begin{eqnarray}
\label{p2_1}\dfrac{\partial(\rho_\sigma{a}_\sigma)}{\partial t}
=\dfrac{\partial}{\partial t}\int\limits_{\mathbb{R}^n}a(t,u)
P(t,x,du)=\int\limits_{\mathbb{R}^n}a(t,u) P_t(t,x,du)=
\nonumber\\
-\int\limits_{\mathbb{R}^n}a(t,u)(a(t,u),\nabla_x
P(t,x,du))+\dfrac12\sigma^2\sum\limits_{k=1}^{n}\dfrac{\partial^2a_\sigma
\rho_\sigma}{\partial x_k^2},
\end{eqnarray}
where $P_t\equiv \frac{\partial }{\partial t} P.$

Further,  we have
\begin{equation}
\begin{array}{l}
\label{p2_2} \dfrac{\partial
(\rho_\sigma\,{a}_{\sigma,k}\,{a}_{\sigma,i})}{\partial x_k}=
{a}_{\sigma,i}(t,x)\dfrac{\partial }{\partial
x_k}\left(\int\limits_{\mathbb{R}^n}\,a_k(t,u)\,P(t,x,du)\right)\,+\\\\
\, \int\limits_{\mathbb{R}^n}\,a_k(t,u)\,P(t,x,du) \,\dfrac{\partial
}{\partial
x_k}\left(\frac{\int\limits_{\mathbb{R}^n}\,a_i(t,u)\,P(t,x,du)}{\int\limits_{\mathbb{R}^n}\,P(t,x,du)}\right)=
\int\limits_{\mathbb{R}^n}\,
a_{\sigma,i}(t,x)\,a_k(t,u)\,P_{x_k}(t,x,du)\,+\nonumber\\\\

\,\int\limits_{\mathbb{R}^n}\,a_k(t,u)\,P(t,x,du)\,
\frac{\int\limits_{\mathbb{R}^n}\,a_i(t,u)\,P_{x_k}(t,x,du)
\,\int\limits_{\mathbb{R}^n}\,P(t,x,du)
\,-\,\int\limits_{\mathbb{R}^n}\,a_i(t,u)\,P(t,x,du)\,\int\limits_{\mathbb{R}^n}\,P_{x_k}(t,x,du)
}{\left(\int\limits_{\mathbb{R}^n}\,P(t,x,du)\right)^2 }=\\\\

\int\limits_{\mathbb{R}^n}\,(a_k(t,u)\, a_{\sigma, i}(t,x)
+a_i(t,u)\, a_{\sigma,k}(t,x)\,-\, a_{\sigma,k}(t,x)\,
a_{\sigma,i}(t,x))\,P_{x_k}(t,x,du),
\end{array}
\end{equation}
$ i,k=1,...,n$, with $P_{x_k} \,\equiv \,\frac{\partial }{\partial
x_k} P.$

Equation \eqref{sist_obw2} follows immediately from \eqref{p2_1} and
\eqref{p2_2}. Thus, Proposition \ref{prop3} is proved.
\bigskip
\begin{corollary}  Before the instant $t_*(u_0)$, the blow up time of the solution to the Cauchy problem
\eqref{burgers_cl},
 the triple $(\bar\rho, \bar u, \bar a)$, which constitutes the limit as $|\sigma|\to 0$
of the triple  $(\rho_\sigma, u_\sigma, a_\sigma)$, solves the
following system:
\begin{equation}\label{pred_rho}\partial_t
\bar\rho+{\rm div}_x(\bar \rho \bar a)=0,\end{equation}
\begin{equation}\label{pred_u}
\partial_t (\bar \rho \bar u)+\nabla_x(\bar\rho \bar u  \bar a)=0,\end{equation}
\begin{equation}\label{pred_a}
\partial_t (\bar\rho \bar a)+\nabla_x(\bar\rho \bar a \otimes \bar a)=0.\end{equation}
\end{corollary}

{\sc Proof.} Equation \eqref{pred_rho} follows from the properties
of parabolic differential equations with a small parameter in front
of the derivatives of second  order (\cite{Freidlin}, Theorem 3.1),
since until the instance $t_*(u_0)$ the coefficients of equation
\eqref{sist_obw1} are differentiable. Equation \eqref{pred_u}
follows from \eqref{pred_rho} and Theorem \ref{T1}, Proposition
\ref{Prop2} implies \eqref{pred_a}.

\begin{remark} System \eqref{pred_rho} and \eqref{pred_a} constitutes the
so called pressureless gas dynamics system , the simplest model
introduced to describe the formation of large structures in the
Universe, see, e.g. \cite{Shand}.
\end{remark}

\begin{remark} As it has been shown in \cite{Roz3} on an example, for discontinuous
solutions to \eqref{burgers_cl} the limits as $\sigma\to 0$ of the
terms $I^a_\sigma$ and $I^u_\sigma$ do not vanish as $\sigma\to 0$
and yield some specific pressure.
\end{remark}

\begin{remark} The method of special stochastic perturbations,
applied here, was used in \cite{Roz1}, \cite{Roz2} for studying
other deterministic problems.
\end{remark}


\end{document}